\documentclass{siamart1116}

\usepackage{color}
\usepackage{graphicx}
\usepackage{amssymb}
\usepackage{amsmath}
\usepackage{subcaption}
\usepackage{multirow}
\usepackage{rotating}
\usepackage{soul}
\usepackage[boxed,vlined,algo2e]{algorithm2e}

\begin{document}

%\begin{frontmatter}

%% Title, authors and addresses

\title{A robust iterative scheme for symmetric  indefinite systems}

\newcommand{\TheTitle}{A robust iterative scheme for symmetric  indefinite systems}
\newcommand{\TheAuthors}{M. Manguo\u{g}lu and V. Mehrmann}

\title{{\TheTitle}\thanks{Submitted DATE.}
\funding{The first author would like to thank Alexander von Humboldt Foundation for the support of a research stay at TU Berlin. The second author thanks Einstein Foundation Berlin for support within the project: \emph{Algorithmic Linear Algebra: High Performance, Numerical Stability, and Error Resiliency.}}}

\headers{\TheTitle}{\TheAuthors}

\author{Murat Manguo\u{g}lu\thanks{Institut f\"ur Mathematik, Technische Universit\"at, Berlin, 10623 Berlin, Germany, and Department of Computer Engineering, Middle East Technical University, 06800 Ankara, Turkey   (\email{manguoglu@ceng.metu.edu.tr})}
\and
Volker Mehrmann\thanks{Institut f\"ur Mathematik, Technische Universit\"at Berlin, 10623 Berlin, Germany  (\email{mehrmann@math.tu-berlin.de})}
}

\maketitle

\begin{abstract}
%% Text of abstract
We propose a two-level nested  preconditioned iterative scheme for solving sparse linear systems of equations in which the coefficient matrix is symmetric and indefinite with relatively small number of negative eigenvalues. The proposed scheme consists of an outer Minimum Residual (MINRES) iteration, preconditioned by an inner Conjugate Gradient (CG) iteration in which CG can be further preconditioned. The robustness of the proposed scheme is illustrated by solving indefinite linear systems that arise in the solution of quadratic eigenvalue problems in the context of model reduction methods for finite element models of disk brakes as well as on other problems that arise in a variety of applications.
\end{abstract}

\begin{keyword}
symmetric indefinite systems, Krylov subspace method,  sparse linear systems,  deflation, preconditioned minimum residual method, preconditioned conjugate gradient method\\
%% keywords here, in the form: keyword \sep keyword
MSC 2000: 65F10, 65F15, 65F50
%% MSC codes here, in the form: \MSC code \sep code
%% or \MSC[2008] code \sep code (2000 is the default)
\end{keyword}

%\end{frontmatter}

%%
%% Start line numbering here if you want
%%
%\linenumbers

%% main text
\section{Introduction}
Symmetric indefinite linear systems
\begin{equation}
\label{eq1}
Ax=b,
\end{equation}
arise in many applications ranging from optimization problems to problems in computational physics, see e.g. \cite{BenGL05,Saa03}.  In this paper we assume that $A\in \mathbb{R}^{n\times n}$ is a sparse, full-rank,  symmetric and indefinite matrix with only few negative eigenvalues. Our motivation to develop a new preconditioned iterative method arises from an application in the automotive industry. In order to control brake squeal, large scale eigenvalue problems are solved via a shift-and-invert Arnoldi method to obtain a reduced model that can be used for parameter studies and optimization, see \cite{GraMQSW16} and Section~\ref{sec:squeal}. We propose the use of a two-level preconditioned iterative method with  a positive definite preconditioner for the solution of the arising linear systems. The basic idea of such a  preconditioned iteration  is well-known. In the context of optimization problems, see \cite{GilMPS92}, a sparse Bunch-Parlett factorization
\begin{equation} \label{pre}
PAP^T  = L D L^T
\end{equation}
is suggested as a solver for the systems involving the indefinite blocks of various preconditioners. Here $P$ is a permutation matrix (with $PP^T = I$), $L$ is a sparse lower triangular matrix (typically  with some fill-in compared to the sparsity pattern of $A$), and $D$ is a block-diagonal matrix that contains either $1\times 1$ or $2\times 2$ blocks. Given such a factorization, one can modify the diagonal matrix $D$ to obtain a positive definite $\tilde{D}$ such that the eigenvalues of $\tilde{D}$ are the absolute values of the eigenvalues of $D$, so that also $M:=L \tilde {D} L^T$ is positive definite. If a diagonal block of $D$ is $1\times 1$ and negative, then one replaces it with its absolute value. Otherwise, it is a symmetric $2\times 2$ block,
\begin{equation}
\begin{bmatrix}
\alpha & \beta \\
\beta & \gamma
\end{bmatrix}.
\end{equation}
and one computes the spectral decomposition
\begin{equation}
\begin{bmatrix}
\alpha & \beta \\
\beta & \gamma
\end{bmatrix} =
\begin{bmatrix}
c & s \\
s & -c
\end{bmatrix}
\begin{bmatrix}
\lambda_1 &  \\
          & \lambda_2
\end{bmatrix}
\begin{bmatrix}
c & s \\
s & -c
\end{bmatrix}
\end{equation}
where $c,s\in \mathbb R$ satisfy $c^2 + s^2 = 1$, and one replaces the $2\times 2$ block with
\begin{equation}
\begin{bmatrix}
\tilde{\alpha} & \tilde{\beta} \\
\tilde{\beta} & \tilde{\gamma}
\end{bmatrix} =
\begin{bmatrix}
c & s \\
s & -c
\end{bmatrix}
\begin{bmatrix}
|\lambda_1| & 0 \\
        0  & |\lambda_2|
\end{bmatrix}
\begin{bmatrix}
c & s \\
s & -c
\end{bmatrix}.
\end{equation}
The matrix $M$, if easily available, is a good  preconditioner for a preconditioned Krylov subspace method, such as the Minimum Residual method (MINRES)~\cite{PaiS75}, since due to the fact that the spectrum of $M^{-1}A$ has only the values $+1, -1$, it would converge in at most $2$ iterations in exact arithmetic  if the factorization is exact. However, this  preconditioner is, in general, not practical for large problems due to fill-in and large storage requirements. In~\cite{GreHL17}, therefore, an incomplete $LDL^T$ factorization based preconditioner for MINRES is proposed.

Another suggestion for a preconditioner of MINRES, proposed in \cite{VerK13},  is the positive definite absolute value of $A$, defined as $|A| := V|\Lambda|V^T$ in which $A=V\Lambda V^T$ is the spectral decomposition of $A$. However, to avoid the high computational complexity of the spectral decomposition, in \cite{VerK13} it is suggested to use  a geometric multigrid  method instead of the absolute value preconditioner and it is illustrated via a model problem that this approach is very effective when the system matrix arises from elliptic partial differential equations.

In our motivating problem, the indefinite matrix arises from  a perturbed wave equation where the resulting linear system depends on parameters and has  the extra property that the number of negative eigenvalues is much smaller than the number of positive eigenvalues. For this class of problems we propose a new two-level iterative scheme that combines the absolute value preconditioner approach with a deflation procedure  and we show that this method is also very effective for a large class of indefinite problems arising in other applications.

\section{A two-level iterative scheme}\label{sec:2level}
In this section we describe a new two-level preconditioned iterative scheme for symmetric indefinite linear systems where the coefficient matrix has only very few negative eigenvalues. The method employs MINRES together with a modified absolute value preconditioner that is constructed via a deflation procedure which,  however, is not carried out explicitly. The linear systems involving the preconditioner are solved again iteratively via the preconditioned Conjugate Gradient (CG)~\cite{FoxHW48} which can be  preconditioned via an incomplete $LU$ ($ILU$) decomposition, see e.g. \cite{Saa03}, of the original coefficient matrix $A$ or any other preconditioner obtained from the original coefficient matrix .  These include but are not limited to Sparse Approximate Inverse~\cite{benzi1996sparse,grote1997parallel}, Algebraic Multigrid~\cite{yang2002boomeramg,brandt1985algebraic}, and banded~\cite{manguoglu2010weighted} preconditioners. We note that just like any iterative method for symmetric linear systems, the proposed scheme would benefit from symmetric permutations and scaling to enhance the numerical stability of the preconditioner as well.  Alternatively, one can use nonsymmetric permutation and scaling~\cite{duff2001algorithms} to further improve the numerical stability of the preconditioner~\cite{benzi2000preconditioning}. However, this would destroy the symmetry and have the disadvantage that robust Krylov subspace methods with short term recurrences are no longer applicable. In this work, we show the robustness of the proposed iterative scheme with classical preconditioners. We illustrate that this MINRES-CG iterative scheme is very effective and more robust than other preconditioned general Krylov subspace methods, such as the restarted Generalized Minimum Residual (GMRES)~\cite{SaaS86}, the stabilized Bi-Conjugate Gradient method (BiCGStab)~\cite{van92},  inner-outer FGMRES-GMRES~\cite{saad1993flexible} using the same (ILU) preconditioner or just modified incomplete $LDL^{T}$ preconditioned MINRES.

As an approximation to the absolute value preconditioner we use
\begin{equation}
M_{mr} := A + 2\hat{V}|\hat{\Lambda}|\hat{V}^T.
\end{equation}
where  $\hat V$ is an approximate invariant subspace of $A$  associated with the (say $k$) negative eigenvalues and $|\hat{\Lambda}|$ is the corresponding absolute value of the diagonal matrix of negative eigenvalues. Since we have assumed that $k$ is much smaller than $n$, the modification (or as it is sometimes called \emph{deflation}) is of small rank.  In each iteration of MINRES applied to~(\ref{eq1}) a system of the form
\begin{equation}
\label{eq2}
M_{mr}z=y
\end{equation}
has to be solved, and again the  preconditioned matrix $M_{mr}^{-1}A$ has only eigenvalues $+1$ or $-1$ so that MINRES with the exact preconditioner converges theoretically again in at most $2$ iterations. However, since $M_{mr}$ is symmetric and positive definite, we propose to use  a preconditioned CG iteration for solving system~(\ref{eq2}) approximately  with an indefinite preconditioner, $M_{cg}$, which is an approximation of the original coefficient matrix itself. Note that the eigenvalues of the preconditioned matrix for CG,  $M^{-1}_{cg} M_{mr}$,  would again be either $+1$ or $-1$ if the exact matrix $A^{-1}$ was used.

Indefinite preconditioning for symmetric positive definite systems is studied in~\cite{freund1991polynomial} where the preconditioned system is solved via a Krylov subspace method other than CG that does not require positive definiteness of the coefficient matrix. Indefinite preconditioning for the CG method is, however, rarely applied with the exception of~\cite{RozS02}, where CG for indefinite systems with indefinite preconditioner is used but it is  assumed that  the preconditioned matrix is positive definite.  In our case, however, this will not be the case.

The first level preconditioner ($M_{mr}$)  is symmetric and positive definite, but dense, so it should not be formed explicitly. On the other hand, the second level preconditioner ($M_{cg}$) is sparse and symmetric but not positive definite. However, the preconditioned CG is still guaranteed {not to break down  (see~\cite[p.~277]{Saa03}}) using an indefinite preconditioner which can be seen as follows. It is well-known, see e.g.~\cite[p.~279]{Saa03}, that preconditioned CG with a preconditioner $M$ applied to a system $Wx=b$ with symmetric positive definite $W$ can be expressed in an indefinite $M$-scalar product by replacing the Euclidean inner products in CG by the $M$-inner products. If $W$ is symmetric positive definite, and $M$ is symmetric indefinite (but invertible), then we can define the indefinite $M$-inner product as $(x,y)_M=(Mx,y)=y^T M x= x^T M y = (y,x)_M$, so $M^{-1}W$ is positive definite with respect to the $M$-inner product, since $(M^{-1}Wx, x)_M > 0$ for all $x\neq 0$.

Given the system $Wx=z$, an initial guess $x_0$, and a preconditioner $M$,  as CG is a projection based Krylov subspace method, the  vectors $x_m$ must satisfy the orthogonality condition
\begin{equation}
\label{orth}
(M^{-1}(z-Wx_m),v)_M = 0 \quad \mbox{for all}\ v \in \hat{\mathcal{K}}_m,
\end{equation}
where $\hat{\mathcal{K}}_m=\mbox{\rm span}\{\hat{r}_0,M^{-1} W\hat{r}_0, ..., (M^{-1}W)^{(m-1)}\hat{r}_0 \}$ and $\hat{r}_0 = M^{-1}r_0$  with $r_0=z-Wx_0$.  Note that~(\ref{orth}) is equivalent to the orthogonality condition of CG without preconditioning
\begin{equation}
(z-Wx_m,v)=0 \quad \mbox{for all}\ v \in \hat{\mathcal{K}}_m
\end{equation}
Therefore, indefinitely preconditioned CG minimizes the error
\begin{equation}
||x_m-x^{*}||_W = \inf_{x\in x_0 + \hat{\mathcal{K}}_m} ||x-x^*||_W,
\end{equation}
in the energy norm defined by the positive definite matrix $W$.

In summary, the proposed scheme consists of three stages. First, an initial preconditioner is obtained from the coefficient matrix $A$, e.g. using the ILU factorization).  Second, we  compute approximations to the  negative eigenvalues and the corresponding invariant subspace. This computation itself may be very expensive even if the invariant subspace has small dimension. However, in our motivating application many linear systems with the same coefficient matrix (or closely related coefficient matrices) need to be solved. Hence, this potentially expensive initial cost is quickly amortized.  This is typical when solving eigenvalue problems  with the shift-and-invert Arnoldi method as in \cite{GraMQSW16}.  The third stage is the iterative solution stage consisting of nested MINRES and CG iterations (Algorithm~\ref{minres_cg}).
Note that while the outer MINRES iterations require matrix-vector multiplications with the original sparse coefficient matrix $A$, the inner CG iterations require matrix-vector multiplications of the form $v = M_{mr} u$  which are efficiently performed by using sparse matrix-vector multiplications and together with dense matrix-vector operations (BLAS Level 2) and vector-vector operations (BLAS Level 1) in the following procedure
\begin{equation}
M_{mr} u= Au + 2(\hat{V}(|\hat{\Lambda}|(\hat{V}^{T}u)))
\end{equation}
The total cost of each such matrix multiplication operation is $O(nnz + kn)$  arithmetic operations where $nnz$, $k$ and $n$  are the number of nonzeros, negative eigenvalues and rows of  $A$, respectively. We note that this extra operation is much more cache friendly than constructing an orthonormal basis in GMRES which relies on inner products (BLAS Level 1). Alternatively, to further speed up the convergence, the proposed scheme can be implemented using RMINRES~\cite{gaul2013framework,wang2007large}, i.e. recycled and deflated MINRES,   as the outer solver instead of plain MINRES. The trade-off would be increased storage and computation requirements due to the necessary orthogonalization against the recycled subspace and the updates of the recycled subspace~\cite{wang2007large}. Furthermore, finding subspaces that lead to improved convergence is considered to be a highly challenging task and application specific~\cite{gaul2013framework}.

\begin{algorithm}
\SetKwFunction{FMain}{MINRES-CG}
\SetKwProg{Pn}{Function}{:}{\KwRet}
\SetKwBlock{Begin}{}{}
\Pn{\FMain{$A$,$b$,$x_0$,$\hat{\Lambda}$,$\hat{V}$}}{
\Begin{
Solve $Ax=b$ via MINRES using the preconditioner $M_{mr}=A + 2\hat{V}|\hat{\Lambda}|\hat{V}^{T}$, the major operations required by each MINRES iteration are:
               \begin{itemize}
                 \item Compute matrix-vector products with $A$.
                 \item Solve $M_{mr} z=y$ via  preconditioned CG using as  preconditioner $M_{cg} = \tilde A$ an approximation of $A$, the major operations required by each CG iteration are:
                 \Begin{ \begin{itemize}
                   \item Compute matrix-vector products: $v=M_{mr}u$
                   \item Solve the system $M_{cg} t = g$
                   \end{itemize}                    }
                 \end{itemize}  }
\KwRet x
}
\caption{Iterative solution stage of MINRES-CG}
\label{minres_cg}
\end{algorithm}

\subsection{Improvement via Sherman-Morrison-Woodbury Formula}
\label{shermanmorrison}
The preconditioner $M_{mr}$ in MINRES-CG is a k-rank update of $A$. Therefore, one can use the Sherman-Morrison-Woodbery formula to express $M^{-1}_{mr}$. Given,
\begin{equation}
M_{mr} := A + 2\hat{V}|\hat{\Lambda}|\hat{V}^T
\end{equation}
after applying the Sherman-Morrison-Woodbury formula and some algebraic manipulations, we obtain,
\begin{equation}
M^{-1}_{mr} := A^{-1} - 2\hat{V}\hat{\Lambda}^{-1}\hat{V}^T.
\end{equation}
Note that since we do not have the exact $A^{-1}$, but use an approximation of it, $M^{-1}_{mr}$ is not positive definite. Still, we can use it as the preconditioner for the $CG$ iterations. In other words, we apply the preconditioner for CG,  $M^{-1}_{cg} = \tilde{A}^{-1} - 2\hat{V}\hat{\Lambda}^{-1}\hat{V}^T$ in which the action of $\tilde{A}^{-1}$ is approximated, such as by an incomplete factorization of $A$. In  the improved scheme, application of the preconditioner involves additional dense matrix-vector multiplications (BLAS Level 2) with a cost of $O(kn)$ arithmetic operations but no additional storage requirement. Hereafter, we refer to this improved version as MINRES-CG$^{*}$.

\section{Application of the two-level method}\label{sec:numerics}
In this section we describe the applications to which we apply the proposed two-level procedure.
\subsection{Finite Element Models of Disk Brakes}\label{sec:squeal}
In the context of noise reduction in disk brakes, reduced order models are determined from  the finite element model  \cite{GraMQSW16} by computing the eigenvalues in the right half plane and close to the imaginary axis of a parametric Quadratic Eigenvalue Problem (QEP)
\begin{equation}
\label{eig.prob}
(\lambda^2 \mathcal{M} + \lambda D_\Omega + K_\Omega )x =0
\end{equation}
in which
\begin{equation}
D_\Omega = D_M + \left(\frac{\Omega_{ref}}{\Omega} -1 \right) D_R + \left(\frac{\Omega}{\Omega_{ref}}\right) D_G
\end{equation}
and
\begin{equation}
K_\Omega = K_E + K_R + \left(\left(\frac{\Omega}{\Omega_{ref}}\right)^2 -1\right)K_g,
\end{equation}
where $\mathcal{M}$ and $K_E$ are symmetric positive definite, $D_G$ is skew-symmetric, $D_M, D_R, K_g$ are symmetric indefinite, and $K_R$ is general~\cite{GraMQSW16}.  Here $\Omega$ denotes the angular velocity of the disk ($2\pi < \Omega < 4\times2\pi$) and $\Omega_{ref}$ is the reference angular velocity.

The QEP is solved by first rewriting it as a linear eigenvalue problem, using a companion linearization of~(\ref{eig.prob}) given by
\begin{equation}
\label{eig.prob.lin}
  \bigg(\begin{bmatrix}
    0 & I \\
    -K_\Omega & -D_\Omega
  \end{bmatrix} - \lambda
   \begin{bmatrix}
   I & 0 \\
   0 & \mathcal{M}
   \end{bmatrix}\bigg)
   \begin{bmatrix}
   x \\
   \lambda x
   \end{bmatrix}
   = 0.
\end{equation}
Audible brake squeal is associated with eigenvalues in the right half plane. For this reason we are interested in those eigenvalues that lie in a rectangular domain in the complex plane given by $-50 < \mbox{\rm Re}(\lambda) < 1,000$ and $-1 < \mbox{\rm Im}(\lambda) < 20,000$ corresponding to the audible range.

Solving the eigenvalue problem~(\ref{eig.prob.lin}) via an eigensolver such as the shift-and-invert Arnoldi method \cite{LehSY98} or Contour Integration based methods \cite{polizzi2009density,sakurai2007cirr}, requires the solution of a shifted linear system of equations in each iteration, see~\cite{GraMQSW16} for details of the eigensolver.  To apply our two-level linear system solver, we consider the solution of the following shifted linear system  with complex shifts ($\gamma$  inside the rectangular domain of interest),
\begin{equation}
\label{complex}
C (x+iy) = f + ig
\end{equation}
where $i=\sqrt{-1}$,  $C = \gamma B - A$, and
\begin{equation}
   B = \begin{bmatrix}
   I & 0 \\
   0 & \mathcal{M}
  \end{bmatrix} ,\ A =
 \begin{bmatrix}
    0 & I \\
    -K_\Omega & -D_\Omega
  \end{bmatrix}.
\end{equation}
In \cite{GraMQSW16} this complex linear system is solved with a sparse complex direct solver. To solve the problem iteratively, we follow~\cite{AxeNA14} and map the complex system~(\ref{complex}) to an equivalent double-size real system.

Splitting into real and imaginary parts $C = \hat{A} + i \hat{B}$ and $\gamma = \gamma_r + i \gamma_i$ with
$\gamma_r = \mbox{\rm Re}(\gamma)$ and $\gamma_i=\mbox{\rm Im}(\gamma)$,
we obtain
\begin{equation}
\hat{A} =
 \begin{bmatrix}
    \gamma_r I & -I \\
    K_\Omega & \gamma_r \mathcal{M}+D_\Omega
  \end{bmatrix},\ \hat{B} =
 \begin{bmatrix}
   \gamma_i I & 0 \\
   0 & \gamma_i \mathcal{M}
  \end{bmatrix},
\end{equation}
for the real and complex parts of $C$, respectively. This leads to the real system
\begin{equation}
\label{blockrealsystem}
\begin{bmatrix}
\hat{B} & -\hat{A} \\
\hat{A} & \hat{B}
\end{bmatrix}
\begin{bmatrix}
x \\
-y
\end{bmatrix}
=
\begin{bmatrix}
g \\
f
\end{bmatrix}.
\end{equation}
which we then solve via a preconditioned Krylov subspace method with preconditioner
\begin{equation}
\label{outer_prec}
M =\begin{bmatrix}
\tilde{B} & -\tilde{A} \\
\tilde{A} & \tilde{B}
\end{bmatrix}
\end{equation}
where
\begin{equation}
\tilde{A} = \begin{bmatrix}
\gamma_r I & -I \\
K_E & \gamma_r \mathcal{M}
\end{bmatrix}
\end{equation}
and $\tilde{B} = \hat{B}$. Note that both $\mathcal{M}$ and $K_E$ are symmetric and positive definite. The preconditioner can be block $LU$ factorized as
\begin{equation}
\begin{bmatrix}
\tilde{B} & -\tilde{A} \\
\tilde{A} & \tilde{B}
\end{bmatrix} =
\begin{bmatrix}
\tilde{B} & 0 \\
\tilde{A} & \tilde{B} + \tilde{A}\tilde{B}^{-1}\tilde{A}
\end{bmatrix}
\begin{bmatrix}
I & -\tilde{B}^{-1}\tilde{A} \\
0 & I
\end{bmatrix}.
\end{equation}
Hence, the major cost in solving systems involving the preconditioner $M$ is the solution of two linear systems where the coefficient matrix is (i) $\tilde{B}$ and (ii) $ S=(\tilde{B}+ \tilde{A}\tilde{B}^{-1}\tilde{A})$, namely the Schur complement. Since the solution of (i) is quite trivial, we only discuss how to solve systems involving the Schur complement matrix, which typically is dense see~\cite{BenGL05}, but in our case it has the factorization
\begin{equation}
\underbrace{\tilde{B}+ \tilde{A}\tilde{B}^{-1} \tilde{A}}_\text{$S$} = \underbrace{\begin{bmatrix}
\mathcal{M}^{-1} & 0 \\
0 & I
\end{bmatrix}}_\text{$S_1$}
\underbrace{\begin{bmatrix}
\left(\gamma_i + \frac{\gamma^2_r}{\gamma_i}\right)\mathcal{M}-\frac{1}{\gamma_i}K_E & -2 \frac{\gamma_r}{\gamma_i}\mathcal{M} \\
2 \frac{\gamma_r}{\gamma_i} K_E & \left(\gamma_i + \frac{\gamma^2_r}{\gamma_i}\right)\mathcal{M}-\frac{1}{\gamma_i}K_E
\end{bmatrix}}_\text{$S_2$}.
\end{equation}
Solving systems involving the Schur complement matrix, therefore, requires two steps: (i) scaling the right hand side vector with $S^{-1}_1$  and (ii) solving systems where the coefficient matrix is $S_2$. Step (i) is again trivial, hence we now look into (ii) which we solve iteratively using a Krylov subspace method where the preconditioner is
\begin{equation}
\tilde{S}_2= \begin{bmatrix}
\left(\gamma_i + \frac{\gamma^2_r}{\gamma_i}\right)\mathcal{M}-\frac{1}{\gamma_i}K_E & 0\\
2 \frac{\gamma_r}{\gamma_i} K_E & \left(\gamma_i + \frac{\gamma^2_r}{\gamma_i}\right)\mathcal{M}-\frac{1}{\gamma_i}K_E
\end{bmatrix},
\end{equation}
 since in our case $||\mathcal{M}||_F\ll||K_E||_F$. Hence, the main cost in solving the block triangular systems lies in the solution
 of
\begin{equation}
\left[\left(\gamma_i + \frac{\gamma^2_r}{\gamma_i}\right)\mathcal{M}-\frac{1}{\gamma_i}K_E \right] u = v,
\end{equation}
or after multiplying both sides of the system by $-\gamma_i$ we obtain
\begin{equation}
\label{inner_level}
[K_E-|\gamma|^2 \mathcal{M}] u = -\gamma_i v,
\end{equation}
where $|\gamma|^2 = \gamma^{2}_i + \gamma^{2}_r$. Even though $\mathcal{M}$ and $K_E$ are symmetric and positive definite, there is no guarantee that the symmetric coefficient matrix $K_E-|\gamma|^2 M$  is positive definite. However, system (\ref{inner_level}) is perfectly suitable for the proposed MINRES-CG scheme, since in our application it only has few negative eigenvalues and they need to be computed only once. Furthermore, the preconditioner (\ref{outer_prec}) is completely independent of the parameters $\Omega$ and $\Omega_{ref}$, and the coefficient matrix of inner systems that have to be solved (\ref{inner_level}) are the same for a given $|\gamma|$. This means that a factorization (incomplete or exact) or an approximation for the coefficient matrix $K_E-|\gamma|^2 \mathcal{M}$ can be computed once and re-used for all values of $\gamma$ of the same absolute value and for all corresponding $\Omega$ values.

Numerical experiments for this class of problems are presented in section~\ref{sec:tests}.

\subsection{Other applications}\label{timdavis}
As further applications we consider all symmetric indefinite problems in the SuiteSparse Matrix Collection~\cite{DavH11} of sizes between $n=1,000$ and $n=50,000$ and with at most $100$ negative eigenvalues.  Since this includes $7$ matrices from  the PARSEC group~\cite{chelikowsky1994finite}, we exclude the $3$ smallest matrices from this group. Furthermore, since shifts around  the so-called Fermi level are also of interest in the PARSEC group of matrices, we shift the largest matrix (SiO) by $A-\sigma I$. For $\sigma$, we chose three values ($0.25$, $0.5$ and $0.75$) which approximately correspond to the gaps in the spectrum of $A$. The properties of these $11$ matrices are given in Table~\ref{UFmatrices}.  Note that we include two examples (*) that arise in finite element discretization of structural problems. These are not full eigenvalue problems but just mass matrices; solving these linear systems is useful if eigenvalues in the inverse mass matrix inner product space are computed.
\begin{table}[htbp]
\caption{Matrices from the SuiteSparse Matrix Collection with application domains and properties ($n$ is matrix dimension, $nnz$ is number of nonzeros and $k$ is number of negative eigenvalues).}
\begin{tabular}{l|r r r l}
Matrix  & $n$ & $nnz$ & $k$ & Application   \\ \hline
Bcsstm$10^{*}$  & $1,086$ & $22,092$ & $54$ & Structural Engineering \\
Bcsstm$27^{*}$ & $1,224$ & $56,126$ & $31$ & Structural Engineering  \\
Nasa$1824$ & $1,824$ & $39,208$ & $20$ & Structural Engineering   \\
%\st{SiH4} & \st{5,041} & \st{171,903} & \st{4} & \st{Real-space pseudopotential method}\\
%\st{SiNa} & \st{5,743} & \st{198,787} & \st{5} & %\st{Real-space pseudopotential method}\\
%\st{Na5} & \st{5,832} & \st{305,630} & \st{4} & \st{Real-space pseudopotential method}\\
Meg$4$ & $5,860$ & $25,258$ & $54$ & RAM Simulation   \\
Benzene & $8,219$ & $242,669$ & $2$ & Real-space pseudopotential method\\
Si$10$H$16$ & $17,077$ & $875,923$ & $41$ & Real-space pseudopotential method\\
Si$5$H$12$ & $19,898$ & $738,598$ & $6$ & Real-space pseudopotential method\\
SiO  & $33,401$ & $1,317,655$ & $8$ & Real-space pseudopotential method\\
SiO($\sigma=0.25$)  & $33,401$ & $1,317,655$ & $16$ & Real-space pseudopotential method\\
SiO($\sigma=0.5$)  & $33,401$ & $1,317,655$ & $26$ & Real-space pseudopotential method\\
SiO($\sigma=0.75$)  & $33,401$ & $1,317,655$ & $41$ & Real-space pseudopotential method\\
\end{tabular}
\label{UFmatrices}
\end{table}

\section{Numerical results}\label{sec:tests}

In this section, we study the robustness of the proposed two-level scheme for indefinite linear systems described in the previous section.  All experiments are performed using MATLAB R201{8}a.

In MINRES-CG, we use the MATLAB eigs function to compute negative eigenvalues and the corresponding eigenvectors. We use an indefinite preconditioner ($M_{cg}$) obtained either by an incomplete $LDL^{T}$ or $LU$  factorization of the coefficient matrix for inner CG iterations. Former is the only suitable preconditioner available in MATLAB which we refer to as $ILU$.  Hence, even though it does not exploit symmetry,  we use it to show the robustness of the proposed scheme in Section~\ref{ilu_uf_comparisons}. For the latter, on the other hand, we use symm-ildl  which is an external package~\cite{GreHL17} that has an interface for MATLAB and is robust. Hereafter, we refer to this preconditioner as $ILDL^{T}$.  We use it to show that the proposed scheme is competitive against other solvers in terms of number of iterations even when a much  more robust preconditioner is used in Sections~\ref{disk_brake_example} and  \ref{ildlt_uf_comparison}.
For a fair comparison, exactly the same preconditioner is used for BiCGStab, GMRES($m$) ($m=20,40,60$  and $120$) as well as another outer-inner scheme with Flexible GMRES (FGMRES) as the outer solver and GMRES as the inner solver~\cite{saad1993flexible}. In Sections~\ref{ildlt_uf_comparison} and \ref{disk_brake_example},  we also use  a MINRES preconditioner with the modified  $ILDL^{T}$ factorization. After computing the $ILDL^T$ factorization, the $D$ matrix is modified as described in~\cite{GilMPS92} in order to obtain a positive definite preconditioner to be used with MINRES.  For FGMRES-GMRES we use a restart value of $120$ for both inner and outer iterations. Iterative solvers, except FGMRES, are the implementations that are available in MATLAB. We note that MATLAB's BiCGStab implementation terminates early before completing a full iteration if the relative residual is already small enough. This counts as a half iteration.  We modified GMRES to stop the iteration based on the true relative residual rather than the preconditioned relative residual.  Storage requirements for MINRES, CG, BiCGStab, FGMRES and GMRES  are given in~\cite{choi2011minres,saad1993flexible,Saa03,SaaS86,van92}, respectively.  In Table~\ref{storage}, we illustrate the storage requirement of each of the iterative solvers via the number of vectors in addition to the coefficient matrix, the preconditioner (i.e. incomplete factors) and  the right hand side vector which are common for all solvers.
\begin{table}[h]
\centering
{
\caption{Total additional memory requirements (number of vectors) of various iterative solver  (not counting $A$ , $M$ and $b$) where $m$ is the restart and $k$ is the number of negative eigenvectors.}
\begin{tabular}{c | c | c | c | c }
MINRES  & MINRES-CG & GMRES & FGMRES-GMRES & BiCGStab \\ \hline
 $7$    &  $11+k$         & $m+2$ &  $3m+4$  &  $6$    \\
\end{tabular}
\label{storage}
}
\end{table}

\subsection{Disk brake example}
\label{disk_brake_example}
In the following we solve (\ref{inner_level}) for the small and large test problems of \cite{GraMQSW16}  of sizes $n=4,669$ and $n=842,638$, respectively, with  $\Omega_{ref}=5$. Note again that (\ref{inner_level}) is independent of $\Omega$. For the first set of experiments we fix the shift $\gamma$ to be the largest value in the range of values of interest, namely $1,000+20,000j$. This also happens to be the most challenging case since the number of negative eigenvalues is also the largest, with
$k=18$ and $k=60$, respectively.

For the proposed scheme, an $ILDL^{T}$ factorization of the coefficient matrix is used as  the preconditioner ($M_{cg}$) of the inner CG iteration. We use the same preconditioner for BiCGStab, GMRES($m$), FGMRES-GMRES  and the modified  $ILDL^{T}$ for MINRES.   For the smaller problem, we also use the $ILU$ factorization with no fill-in (i.e. $ILU(0)$) preconditioner of MATLAB.

For all experiments a moderate outer stopping tolerance of relative residual norm less than or equal to $10^{-3}$ is used. For the MINRES-CG  and FGMRES-GMRES schemes the inner stopping tolerance is $10^{-2}$. For all methods, the maximum (total) number of iterations are $2,000$ and $15,000$ for small and large problems, respectively. In all experiments, the right hand side vector is a random vector of size $n$.

The required number of iterations for the proposed scheme as well as for baseline algorithms are given in Table~\ref{table1} for solving the small problem using $ILU(0)$, $ILDL^{T}(1,10^{-2})$ and $ILDL^{T}(1,10^{-3})$ preconditioners. GMRES($20$) reaches the maximum number of iterations without converging ($\dagger$) irrespective of the preconditioner. When the preconditioner is $ILU(0)$, BiCGStab converges but it requires twice as many iterations as MINRES-CG, while all other  solvers reach the maximum number of iterations without converging. Using $ILDL^{T}(1,10^{-2})$ and $ILDL^{T}(1,10^{-3})$ as the preconditioners, BiCGStab , MINRES and GMRES($m$) converge for  $m=120$ and $m=40,60,120$, respectively.

In Table~\ref{table2}, results are presented for solving the large problem using the seven iterative methods with the preconditioners  $ILDL^{T}(4,10^{-4})$, $ILDL^{T}(5,10^{-5})$ as well as  $ILDL^{T}(5,10^{-6})$. Note that a much smaller dropping tolerance is required for the large problem. Incomplete factors contain $117.1$ , $144.5$, and $146.8$  nonzeros per row, respectively, which is relatively small considering that the complete $LDL^{T}$ factorization would produce $558.4$ nonzeros per row. In fact, incomplete factorization may not be an efficient preconditioner for this problem. However, we still include these results here only to show the robustness of the proposed scheme in terms of number of iterations. While for all preconditioners GMRES($m$), FGMRES-GMRES, BiCGStab and MINRES reach the maximum number of iterations without converging, MINRES-CG still converges in $4$ outer iterations albeit with a large number of inner iterations.

\begin{table}[htbp]
\centering
\caption{{ Required number of iterations using various preconditioners and iterative methods for the small system ($\dagger$: maximum number of iterations is reached without convergence)}}
\begin{tabular}{l|l c c c}
Preconditioner & Solver & Outer its.   & Inner its. (Avg.) & Total its. \\ \hline
\multirow{7}{*}{$ILU(0)$}& BiCGStab  & $1,421.5$  & -  & $1,421.5$  \\
& GMRES($20$)  & $\dagger$  & - & $\dagger$  \\
& GMRES($40$)  & $\dagger$  & - & $\dagger$ \\
& GMRES($60$)  & $\dagger$  & - & $\dagger$ \\
& GMRES($120$)  & $\dagger$ & - & $\dagger$  \\
& FGMRES-GMRES & $\dagger$& $\dagger$ & $\dagger$   \\
& MINRES-CG  & $4$ & $177.75$  & $711$ \\ \hline
\multirow{8}{*}{$ILDL^{T}(1,10^{-2})$}
& BiCGStab  & $1,337.5$  & -  & $1,337.5$  \\
& GMRES($20$)  & $\dagger$ & - & $\dagger$  \\
& GMRES($40$)  & $\dagger$  & - & $\dagger$ \\
& GMRES($60$)  & $\dagger$  & - & $\dagger$ \\
& GMRES($120$)  & $7(103)$ & - &  $823$  \\
& MINRES  & $591$  & - &     $591$ \\
& FGMRES-GMRES & $1(3)$   & $259.3$ & $778$  \\
& MINRES-CG  & $4$ & $199.5$  & $798$ \\
\hline
\multirow{8}{*}{$ILDL^{T}(1,10^{-3})$}
& BiCGStab  & $143$  & -  & $143$  \\
& GMRES($20$) & $\dagger$  & - & $\dagger$  \\
& GMRES($40$)  & $8(40)$  & - & $200$ \\
& GMRES($60$)  & $2(7)$  & - & $67$ \\
& GMRES($120$)  & $1(62)$ & - &  $62$ \\
& MINRES  & $164$ & - & $164$ \\
& FGMRES-GMRES & $1(4)$   & $46.3$ & $185$  \\
& MINRES-CG  & $4$ & $32.8$  & $131$ \\
\end{tabular}
\label{table1}
\end{table}

\begin{table}[htbp]
\centering
\caption{{ Required number of iterations using various preconditioners and iterative methods for the large system ($\dagger$:  maximum number of iterations is reached without convergence)}}
\begin{tabular}{l|l c c c}
Preconditioner & Solver & Outer its.   & Inner its. (Avg.) & Total its. \\ \hline
\multirow{8}{*}{$ILDL^{T}(4,10^{-4})$}& BiCGStab  & $\dagger$  & -  & $\dagger$  \\
& GMRES($20$)  & $\dagger$  & - & $\dagger$ \\
& GMRES($40$)  & $\dagger$  & - & $\dagger$ \\
& GMRES($60$)  & $\dagger$  & - & $\dagger$ \\
& GMRES($120$)  & $\dagger$ & - & $\dagger$ \\
& MINRES &  $\dagger$  &  -  & $\dagger$         \\
& FGMRES-GMRES & $\dagger$  & $\dagger$   & $\dagger$   \\
& MINRES-CG  & $4$ & $3,032$  & $12,128$ \\ \hline
\multirow{8}{*}{$ILDL^{T}(5,10^{-5})$}& BiCGStab  & $\dagger$  & -  & $\dagger$  \\
& GMRES($20$)  & $\dagger$  & - & $\dagger$  \\
& GMRES($40$)  & $\dagger$  & - & $\dagger$ \\
& GMRES($60$)  & $\dagger$  & - & $\dagger$ \\
& GMRES($120$) & $\dagger$ & - &  $\dagger$ \\
%& {\color{red} MINRES}&  {\color{red}$15,789$}  &  - &  {\color{red}$15,789$}           \\
& MINRES&  $\dagger$  &  - &  $\dagger$      \\
& FGMRES-GMRES & $\dagger$  & $\dagger$   & $\dagger$   \\
& MINRES-CG  & $4$ & $2,221$  & $8,884$ \\ \hline
\multirow{8}{*}{$ILDL^{T}(5,10^{-6})$}& BiCGStab  & $\dagger$  & -  & $\dagger$  \\
& GMRES($20$)  & $\dagger$  & - & $\dagger$  \\
& GMRES($40$)  & $\dagger$  & - & $\dagger$ \\
& GMRES($60$)  & $\dagger$  & - & $\dagger$ \\
& GMRES($120$) & $\dagger$ & - &  $\dagger$      \\
& MINRES &  $\dagger$  &  -  & $\dagger$   \\
& FGMRES-GMRES & $\dagger$  & $\dagger$   & $\dagger$   \\
& MINRES-CG  & $4$ & $2,242.8$  & $8,971$ \\
\end{tabular}
\label{table2}
\end{table}

\begin{figure}
\includegraphics[scale=0.6]{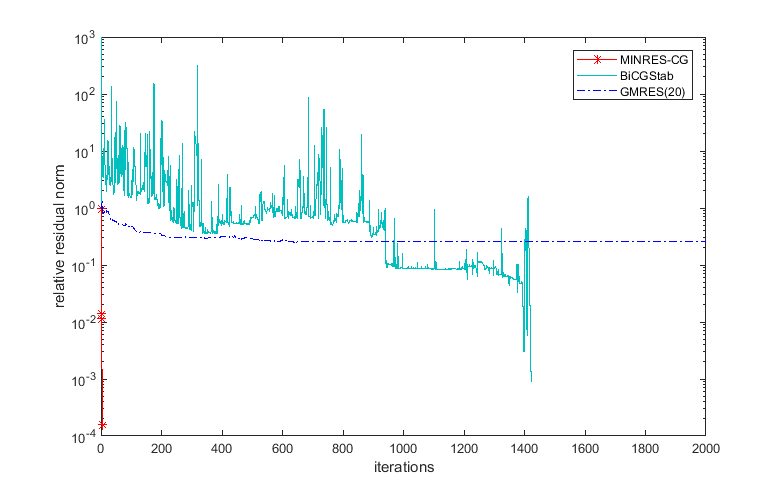}
\caption{The relative residual history for MINRES-CG, BiCGStab and GMRES(20).}
\label{resid_hist_omega5}
\end{figure}

In Figure~\ref{resid_hist_omega5}, the relative residual history is given when the $ILU(0)$ preconditioner is used for three algorithms for the small test problem. Note that for MINRES-CG the relative residual is only available at each outer iteration. Hence, only those are presented in the figure.

\begin{figure}
\begin{subfigure}{.5\textwidth}
  \centering
  \includegraphics[scale=0.44]{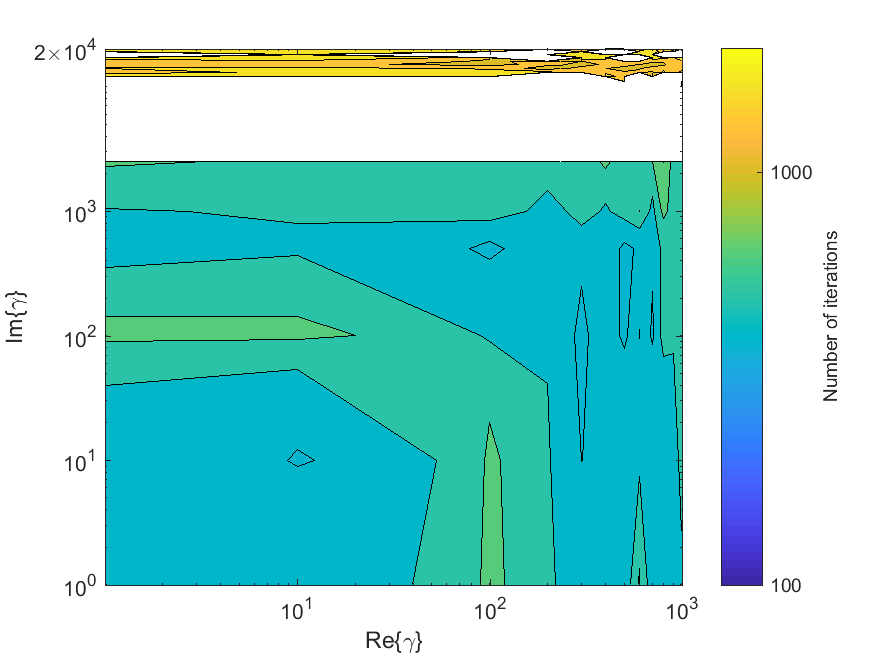}
  \caption{BiCGStab}
%  \label{fig:sfig1}
\end{subfigure}%
\begin{subfigure}{0.5\textwidth}
  \centering
  \includegraphics[scale=.44]{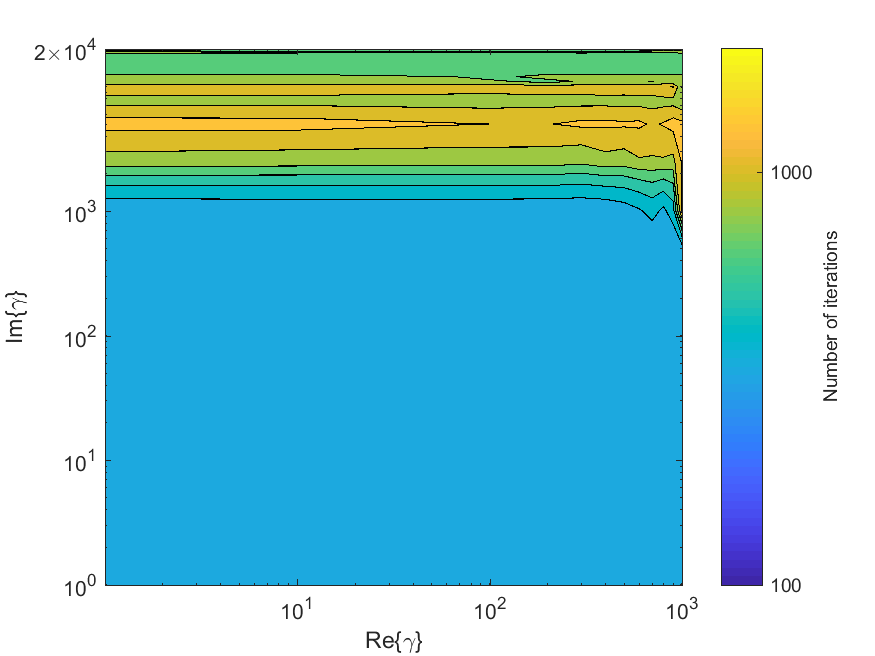}
  \caption{MINRES-CG}
%  \label{fig:total_omegaref5}
\end{subfigure}

\caption{Total number of iterations for BiCGStab and MINRES-CG using the preconditioner $ILU(0)$. White color indicates that the method failed to converge. GMRES($20$) fails for all shifts.}
\label{fig:total_omegaref5}
\end{figure}

\begin{figure}
\begin{subfigure}{.5\textwidth}
  \centering
  \includegraphics[scale=0.44]{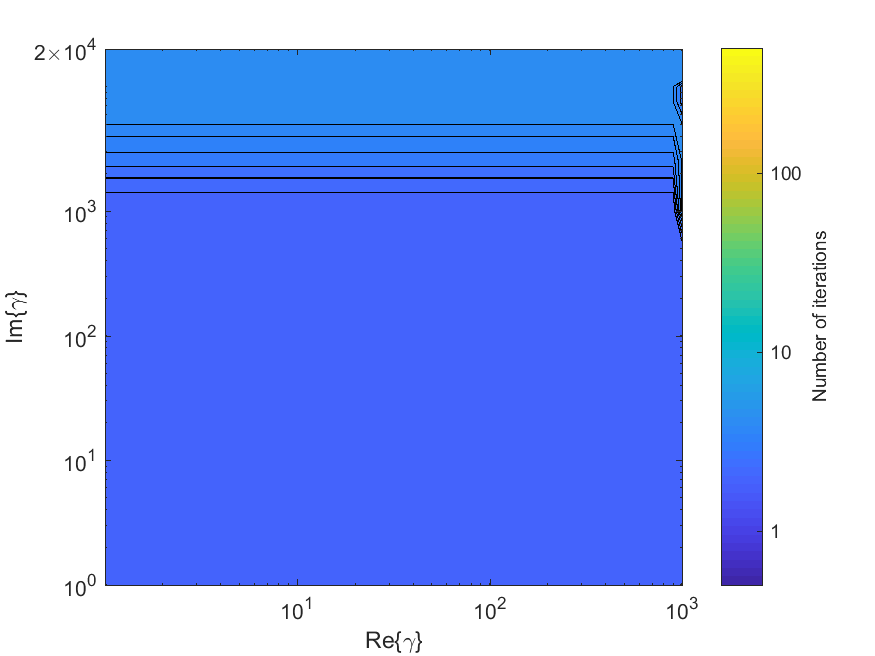}
  \caption{Outer}
%  \label{fig:sfig1}
\end{subfigure}%
\begin{subfigure}{.5\textwidth}
  \centering
  \includegraphics[scale=.44]{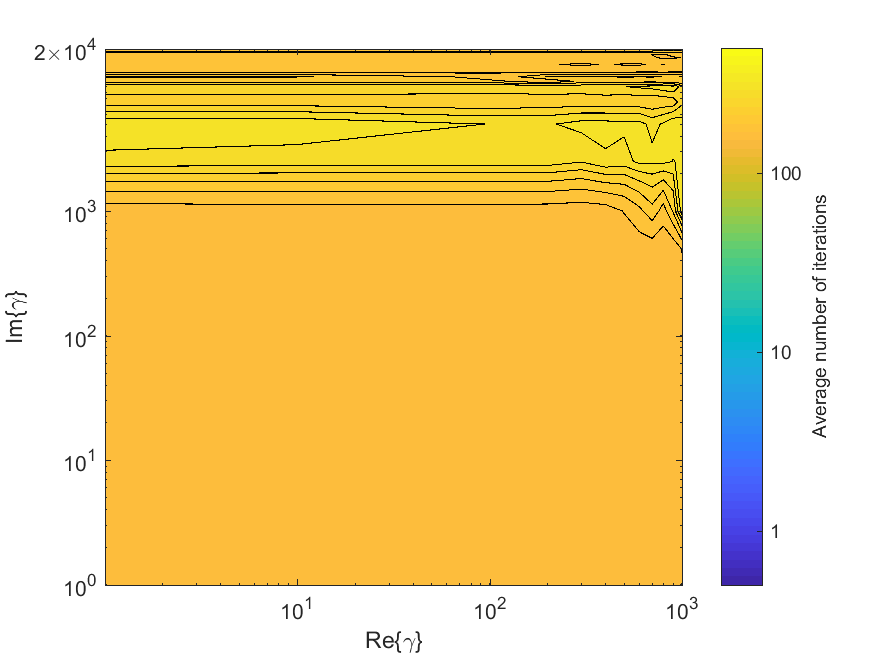}
  \caption{Inner}
%  \label{fig:sfig2}
\end{subfigure}

\caption{Number of outer (MINRES) iterations and average number of inner (CG) iterations for MINRES-CG using the preconditioner $ILU(0)$.}
\label{fig:omegaref5}
\end{figure}

As second application we fix the preconditioner to be $ILU(0)$ and vary the shift $\gamma$ in the complex domain of interest for the small test problem.  Here $\gamma$ is a parameter that we change in the context of the eigenvalue problem. It is of interest to see how  the method behaves as $\gamma$ is varied. In Figure~\ref{fig:total_omegaref5}, the total number of iterations is presented. Preconditioned BiCGStab fails to converge for some values of $\gamma$ (shown as the white area in the figure) while MINRES-CG converges for all  $\gamma$ values.  Figure~\ref{fig:omegaref5} depicts  the number of outer iterations and the average number of inner iterations for MINRES-CG.

\subsection{Test cases from the SuiteSparse matrix collection}
\label{iluprecondsection}
In this subsection, the results are presented for systems that are obtained from the SuiteSparse Matrix Collection. In the first set we compare the proposed method against the classical general iterative schemes GMRES($m$) and BiCGStab as well as a nested FGMRES-GMRES scheme using an incomplete LU factorization based preconditioner. In the second set, we compare the proposed method against the same iterative methods as in the first set but using the $ILDL^{T}$ preconditioner, and MINRES with the modified $ILDL^{T}$ preconditioner.

\subsubsection{Comparison against $ILU$ preconditioner}
\label{ilu_uf_comparisons}

We use $ILU(0)$ for all cases except Meg$4$ where incomplete LU factorization fails due to a zero pivot. Therefore, we use the modified incomplete LU factorization in MATLAB with $10^{-2}$ dropping tolerance (i.e. $MILU(10^{-2})$) for this case only. Since in practice GMRES is  always used with a value for the restart ($m$) we choose  a restart value of $m=20,40,60$ and $120$.   In FGMRES-GMRES, we use a restart of $120$ for both inner and outer iterations.  We stop the iterations when the relative residual norm is less than $10^{-5}$ for all cases. The inner iteration stops when the relative residual norm is less than $10^{-3}$ for CG and GMRES, in MINRES-CG and FGMRES-GMRES, respectively.  Both in MINRES-CG and BiCGStab  iterations stop when the true relative residual is less than the tolerance. For preconditioned GMRES the available residual is only the preconditioned residual. In order to have a fair comparison, we explicitly compute the true residual at each GMRES iteration and stop the iteration based on the true relative residual norm. For all methods, the maximum (total) number of iterations is $20,000$.

In Table~\ref{detailed_iters}, the detailed number of iterations for $ILU$ preconditioned MINRES-CG, GMRES($m$), FGMRES-GMRES and BiCGStab are given. GMRES($20$) fails in $6$ cases out of $11$. For $bcsstm10$, GMRES($20$) stagnates ($\ddagger$), while for $5$ other cases (namely $bcsstm27$, $nasa1824$, $Si10H16$, $SiO(\sigma=0.25)$ and $Sio(\sigma=0.75)$), the maximum number of iterations is reached without convergence ($\dagger$). If the restart is increased to $40$, $60$ and $120$, GMRES($m$), fails in $4$, $3$ and $2$ cases, respectively. BiCGStab fails for $bcsstm27$ and $Meg4$ due to the maximum of iterations being reached without convergence ($\dagger$) and a scalar quantity became too large or too small during the iteration ($\ast$), respectively. FGMRES-GMRES fails in $3$ cases due to the maximum of iterations being reached without convergence ($\dagger$).   The proposed MINRES-CG method does not fail in any of the test problems. Although the cost per iteration is different for each method, the total number of iterations are presented in Table~\ref{total_iters}. For the cases they do not fail, GMRES($120$) and FGMRES-GMRES require fewer number of iterations than MINRES-CG but they also require more storage.  In $4$ cases MINRES-CG requires fewer iterations than BiCGStab. It is possible to improve the total number of iterations of MINRES-CG via using the algorithm described in Section~\ref{shermanmorrison}, Table~\ref{new_vs_old} shows the improved number of iterations which is significant especially for the cases where the inner or outer number of iterations are high.

\begin{figure} [h]
%\begin{subfigure}[b]{0.99\textwidth}
\centering
  \includegraphics[height=3.2in]{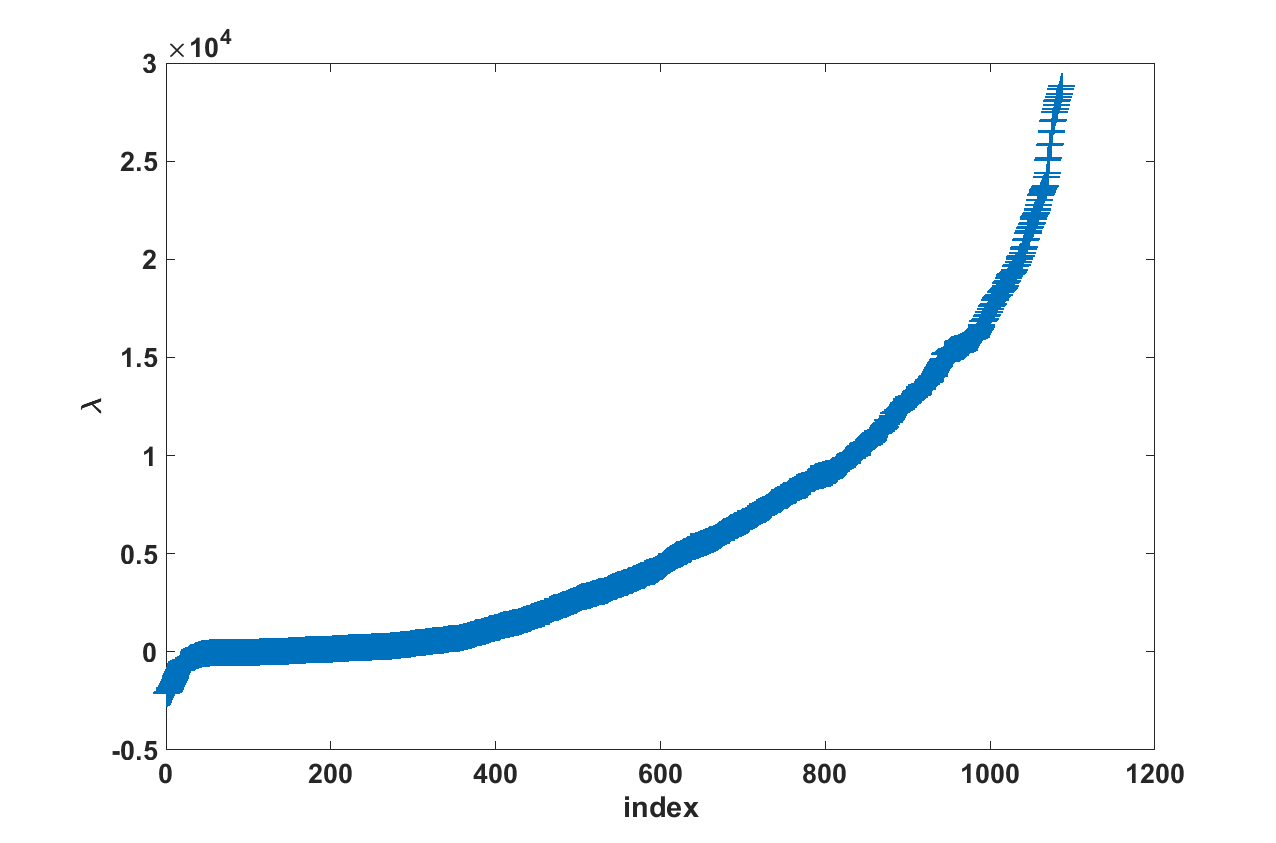}
  \caption{Eigenvalues of A ($bccstm10$)}
  \label{fig:eigs_a}
%\end{subfigure}% \
\end{figure}

\begin{figure}[h]
%\begin{subfigure}[b]{0.99\textwidth}
\centering
  \includegraphics[height=3.4in]{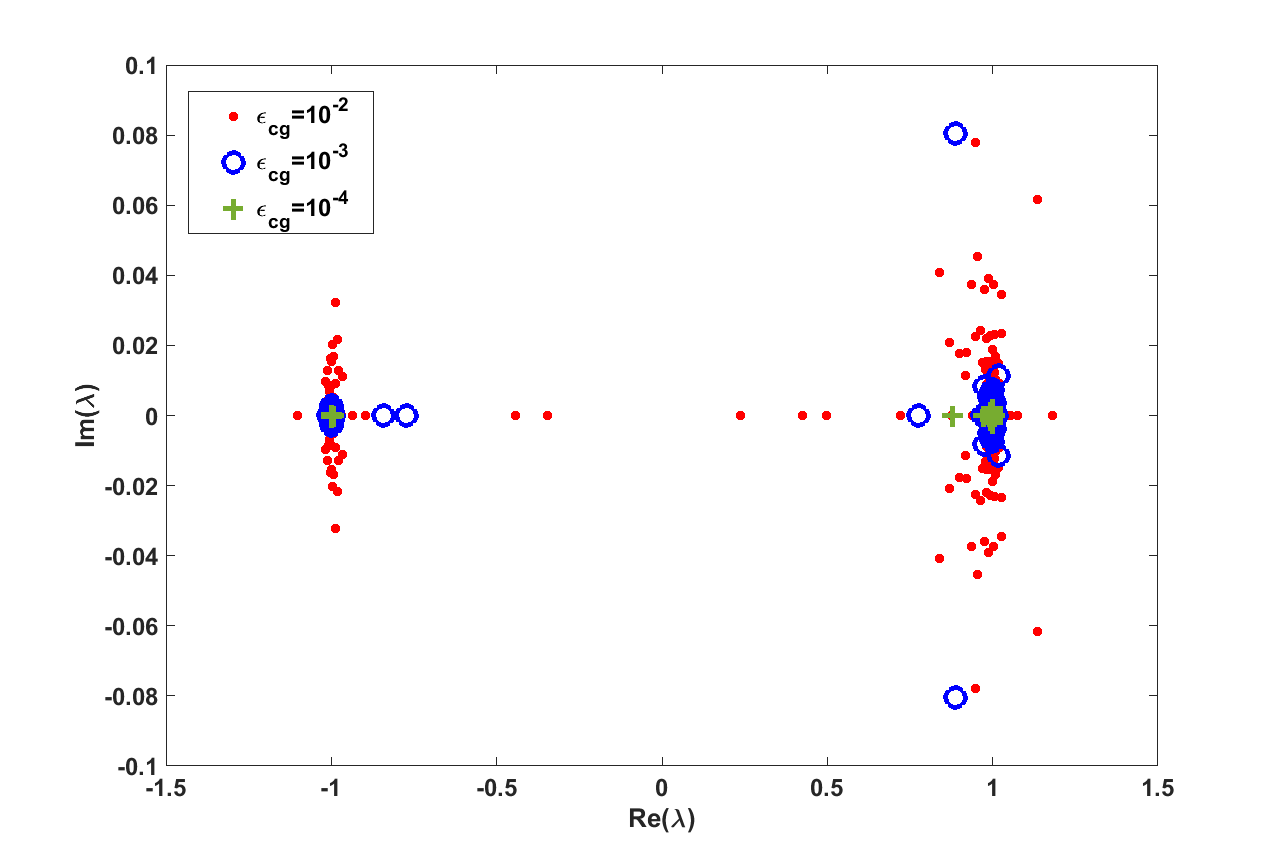}
  \caption{Eigenvalues of $M_{mr}^{-1}A$ ($\epsilon_{cg}=10^{-2}$, $10^{-3}$ and $10^{-4}$) ($bccstm10$)}
  \label{fig:eigs_c}
%\end{subfigure}
\end{figure}

In order to study the effect of the inner stopping tolerance on the eigenvalues of the preconditioned matrix, we explicitly compute $M_{mr}^{-1}A$ using preconditioend CG iterations with stopping tolerances of $10^{-2}$, $10^{-3}$ and $10^{-4}$.  In Figure~\ref{fig:eigs_c},  a clear clustering of  eigenvalues of the preconditioned matrix $M_{mr}^{-1}A$ is visible around $+1$ and $-1$ for $bcsstm10$, while the unpreconditioned coefficient matrix had no clustering of eigenvalues (see Figure~\ref{fig:eigs_a}).  As expected, the clustering around $-1$ and $+1$ improves as the stopping tolerance for the inner CG iterations is decreased.

\begin{table}[h]
\centering
\caption{Comparison of MINRES-CG and MINRES-CG$^{*}$, the improvement via  the Sherman-Morrison-Woodbury formula  is given in the second column, both are using the same $ILU$ preconditioner.}
\begin{tabular}{l | c c | c c }
& \multicolumn{2}{|c|}{MINRES-CG} & \multicolumn{2}{|c}{ {MINRES-CG$^{*}$}}  \\
Name & MINRES  & CG (Avg.) & {MINRES}  & {CG (Avg.)}  \\ \hline
Bcsstm$10$ & $4$ & $650.5$ & ${4}$ &  ${470}$   \\
Bcsstm$27$ & $5$ & $3,186.4$ & ${4}$ & ${1,339.5}$   \\
Nasa$1824$ & $4$ & $455.3$ & ${4}$ & ${309.3}$ \\
%iH$4$ & $4$ & $31.3$ & ${4}$ & ${34}$  \\
%iNa & $4$ & $26.3$ & ${4}$ & ${25.3}$  \\
%a$5$ & $4$ & $28.5$ & ${3}$ & ${22.7}$  \\
Meg$4$ & $16$ & $18.6$ & ${4}$ & ${1.5}$  \\
Benzene & $3$ & $24.3$ & ${3}$ &${22.7}$  \\
Si$10$H$16$ & $4$ & $831$ & ${4}$ & ${893.5}$  \\
Si$5$H$12$ & $4$ & $53.8$ & ${4}$ & ${47.3}$  \\
SiO  & $4$ & $50.5$ & ${4}$ &  ${48.8}$  \\
SiO($\sigma=0.25$) & $4$  & $259$   & ${4}$ & ${141.8}$ \\
SiO($\sigma=0.5$) & $4$  & $94.3$   & ${4}$ & ${80}$ \\
SiO($\sigma=0.75$) & $4$  & $179.5$   & ${4}$ & ${184}$
\end{tabular}
\label{new_vs_old}
\end{table}

\newpage

\begin{sidewaystable}[htbp]
\centering
\caption{Number of iterations  using $ILU$}
\begin{tabular}{l |c r | c c c  c | c c | c }
& \multicolumn{2}{|c|}{MINRES-CG} & \multicolumn{4}{|c|}{GMRES($m$)} & \multicolumn{2}{|c|}{{FGMRES($m_1$)-GMRES($m_2$)}} & BiCGStab  \\
Name & MINRES  & CG & $m=20$ & {$m=40$} & {$m=60$} & {$m=120$} & {$m_1=120$} & {$m_2=120$} &  \\ \hline
Bcsstm$10$ & $4$ & $650.5$ & $\ddagger$ &  {$\ddagger$} &   {$12(60)$} & {$1(80)$} & {$1(57)$} & {$14.5$} & $443.5$  \\
Bcsstm$27$ & $5$ & $3,186.4$ & $\dagger$ & {$\dagger$}  &  {$\dagger$} & {$\dagger$} & {$\dagger$}  & {$\dagger$} &  $\dagger$ \\
Nasa$1824$ & $4$ & $455.3$ & $\dagger$ & {$101(14)$} & {$15(30)$} & {$4(78)$} & {$1(6)$} & {$92$} & $483$ \\
%SiH$4$ & $4$ & $31.3$ & $5(2)$ & {$2(16)$} & {$1(45)$} & {$1(45)$} & {$1(3)$} & {$31.3$} & $57$ \\
%SiNa & $4$ & $26.3$ & $7(1)$ & {$2(2)$} & {$1(42)$} & {$1(42)$}  &  {$1(2)$} & {$32.0$} & $50.5$ \\
%Na$5$ & $4$ & $28.5$ & $12(1)$ & {$2(21)$} & {$1(51)$} & {$1(51)$} & {$1(2)$} & {$38.0$} & $45.5$ \\
Meg$4$ & $16$ & $18.6$ & $5(20)$ & {$1(38)$}  & {$1(38)$} & {$1(38)$} & {$1(3)$} & {$28.3$} & $\ast$ \\
Benzene & $3$ & $24.3$ & $6(4)$ & {$2(1)$} & {$1(41)$} & {$1(41)$} & {$1(2)$} & {$30.5$}  & $37.5$ \\
Si$10$H$16$ & $4$ & $831$ & $\dagger$ & {$\dagger$} & {$\dagger$} & {$\dagger$} & {$\dagger$} & {$\dagger$} & $6,956.5$ \\
Si$5$H$12$ & $4$ & $53.8$ & $347(19)$ & {$9(18)$} & {$3(14)$}  & {$1(80)$} & {$1(3)$} & {$66.7$} &  $108.5$ \\
SiO  & $4$ & $50.5$ & $106(5)$ &  {$6(15)$} & {$2(18)$} & {$1(70)$} & {$1(2)$} & {$55$} & $114$ \\
{SiO ($\sigma=0.25$)} & {$4$} & {$259$} & {$\dagger$} & {$\dagger$} & {$\dagger$} & {$25(120)$} & {$\dagger$} & {$\dagger$} & {$2,109.5$}  \\
{SiO ($\sigma=0.5$)} & {$4$} & {$94.3$} & {$492(20)$} & {$48(31)$} & {$12(59)$}  & {$2(111)$} &  {$1(3)$} & {$116.3$} & {$621.5$}   \\
{SiO ($\sigma=0.75$)} & {$4$} & {$179.5$} & {$\dagger$} &  {$51(33)$}  & {$20(24)$} & {$5(56)$}  & {$1(3)$} & {$232.7$} & {$1,007$}  \\
\end{tabular}
\label{detailed_iters}
\end{sidewaystable}

\begin{sidewaystable}[htbp]
\centering
\caption{Total number of iterations using $ILU$}
%\footnotesize
\begin{tabular}{l |c| c c c c | c | c }
&MINRES-CG&\multicolumn{4}{|c|}{GMRES($m$)}&FGMRES($120$)&BiCGStab\\
Name&& $20$ &$40$&$60$& $120$&-GMRES($120$)&\\\hline
%\begin{table}
%\caption{Total number of iterations using $ILU$}
%\begin{tabular}{l| r r r}
%Name  & MINRES-CG  & GMRES($20$) & BiCGStab  \\ \hline
Bcsstm$10$  & $2,642$ & $\ddagger$ & $\ddagger$ & $720$ & $80$ & $826$  & $443.5$ \\
Bcsstm$27$ & $15,932$ &  $\dagger$ & $\dagger$ &  $\dagger$ & $\dagger$ & $\dagger$ &  $\dagger$ \\
Nasa$1824$ & $1,821$ & $\dagger$ & $4,014$  & $870$ & $438$ & $552$ & $483.0$ \\
%iH$4$ & $125$ & $82$ & $56$ & $45$ & $45$ & $94$ & $57.0$ \\
%iNa & $105$ & $121$ & $56$ & $42$ & $42$ & $64$ & $50.5$ \\
%a$5$ & $114$ & $221$ & $61$ & $51$ & $51$ & $76$ &   $45.5$ \\
Meg$4$ & $297$ & $100$ & $38$ & $38$ & $38$ & $85$ &  $\ast$ \\
Benzene & $73$ & $104$ & $41$ & $41$ & $41$ & $61$ &  $37.5$ \\
Si$10$H$16$ & $3,324$ & $\dagger$ & $\dagger$  & $\dagger$ &  $\dagger$ & $\dagger$ &   $6,956.5$ \\
Si$5$H$12$ & $215$ & $6,939$ & $338$ & $134$ & $80$ & $200$ &$108.5$ \\
SiO  & $202$ & $2,105$ & $215$ & $78$ & $80$ & $110$ & $114$ \\
SiO ($\sigma=0.25$) & $1,036$ & $\dagger$ & $\dagger$ & $\dagger$ & $\dagger$ & $\dagger$ & $2,109.5$ \\
SiO ($\sigma=0.5$) & $377$ & $9,840$ & $1,911$ & $719$ & $231$  &  $349$ & $621.5$ \\
SiO ($\sigma=0.75$) & $718$ & $\dagger$ & $2,033$ & $1,164$ & $536$ & $689$ & $1,007$
\end{tabular}
\label{total_iters}
\end{sidewaystable}

\subsubsection{Comparisons using incomplete $LDL^{T}$}
\label{ildlt_uf_comparison}
In this subsection, we compare the proposed scheme against a robust incomplete $LDL^T$ (Bunch-Parlett) factorization~\cite{Fis11}.

We use the $ILDL^T$ implementation of~\cite{GreHL17} in MATLAB which  computes the incomplete Bunch-Parlett factorization of the coefficient matrix. The default parameters are $3$ and $10^{-3}$ for the level of fill-in and the dropping tolerance, respectively. Furthermore, it uses the Approximate Minimum Degree reordering, Rook pivoting and  scaling to improve the numerical stability of the incomplete factors by default. Note that all of those enhancements that are implemented in $ILDL^T$ make the preconditioner much more robust than the $ILU(0)$  preconditioner.  In the following experiments all methods are applied to the permuted and scaled linear systems. We use the modified  $ILDL^T$ factorization as a preconditioner for MINRES. To have a fair comparison, the same $ILDL^T$ factorization (without the modification) is used as the preconditioner for GMRES($m$), FGMRES-GMRES,  BiCGStab and as the inner preconditioner for MINRES-CG. Stopping tolerances and the maximum number of iterations allowed are set exactly the same as in Section~\ref{iluprecondsection}.

In  Table~\ref{detailed_iters_ildl} the total number of iterations for all methods are given. Even though it is a much more robust preconditioner,  MINRES preconditioned with the modified $ILDL^T$ preconditioner stagnates ($\ddagger$) for $bcsstm27$. For the same problem  BiCGStab, FGMRES-GMRES and GMRES($m$) (for all restart values $m=20,40,60,120$) reach the maximum number of iterations without converging ($\dagger$). On the other hand,  MINRES-CG converges in all problems which confirms the robustness of the proposed scheme. In Table~\ref{total_iters_ildl}, the total number of iterations for all methods are given. GMRES($m$) with larger restart values and BiCGStab require the fewest number of iterations. FGMRES-GMRES requires fewer iterations than MINRES-CG. On the other hand, MINRES requires more iterations than MINRES-CG in $4$ cases, and the required number of iterations are marginally better than that of MINRES-CG for $3$ other cases.

\begin{sidewaystable}[htbp]
\centering
\caption{Number of iterations using $ILDL^T$  preconditioner (MINRES uses the modified spd preconditioner)}
\begin{tabular}{l |c r | c c c  c | c c | c | c}
& \multicolumn{2}{|c|}{MINRES-CG} & \multicolumn{4}{|c|}{GMRES($m$)} & \multicolumn{2}{|c|}{FGMRES($m_1$)-GMRES($m_2$)} & BiCGStab  & MINRES \\
Name & MINRES  & CG & $m=20$ & $m=40$ & $m=60$ & $m=120$ & $m_1=120$ & $m_2=120$ &  \\ \hline
Bcsstm$10$ & $4$ & $10.5$ & $1(10)$ & {$1(10)$} & {$1(10)$} & {$1(10)$} & {$1(4)$} & {$5.5$} & $8.5$ & $37$   \\
Bcsstm$27$ & $4$ & $2,127$ & $\dagger$ & {$\dagger$}  & {$\dagger$} & {$\dagger$} & {$\dagger$} & {$\dagger$} & $\dagger$ & $\ddagger$  \\
Nasa$1824$ & $4$ & $19.5$ & $2(11)$ & {$1(23)$} & {$1(23)$}  & {$1(23)$}  & {$1(21)$}  & {$3.5$}   & $32$ & $48$  \\
%iH$4$ & $4$ & $8.3$ & $1(9)$ & {$1(9)$} & {$1(9)$} & {$1(9)$} & {$1(3)$} & {$4.0$} & $5$ & $18$\\
%iNa & $4$ & $10.5$ & $1(9)$ & {$1(9)$} & {$1(9)$} & {$1(9)$} & {$1(3)$} & {$3.6666$} & $4.5$ & $23$ \\
%a$5$ & $4$ & $9.0$ & $1(12)$ & {$1(12)$} & {$1(12)$} & {$1(12)$}  & {$1(3)$} & {$6.0$} & $8$ & $25$  \\
Meg$4$ & $2$ & $1.5$ & $1(1)$ & {$1(1)$}  & {$1(1)$} & {$1(1)$} & {$1(1)$} & {$1$} & $0.5$ & $2$ \\
Benzene & $3$ & $7.7$ & $1(9)$ & {$1(9)$} & {$1(9)$} & {$1(9)$} & {$1(3)$} & {$4.7$} & $4.5$ & $12$ \\
Si$10$H$16$ & $4$ & $72$ & $7(12)$ & {$2(13)$} & {$1(42)$} & {$1(42)$} & {$1(3)$} & {$27$}  & $49.5$ & $1,192$\\
Si$5$H$12$ & $4$ & $17.3$ & $1(17)$ & {$1(17)$} &{$1(17)$} & {$1(17)$} & {$1(3)$} &  {$10.7$} & $11.5$ & $34$\\
SiO  & $4$ & $18$ & $1(17)$ & {$1(17)$} & {$1(17)$} & {$1(17)$} & {$1(3)$} & {$10.7$} & $11.5$ & $51$ \\
{SiO($\sigma=0.25$)} & {$4$}  & {$24.8$}  & {$4(3)$} & {$1(31)$}  & {$1(31)$}  & {$1(31)$} & {$1(3)$} & {$18$}  & {$40.5$}  & {$80$} \\
{SiO($\sigma=0.5$)} & {$4$}  & {$43.3$}  & {$12(20)$} & {$2(2)$}  & {$1(41)$}  & {$1(41)$} & {$1(3)$} & {$28.3$}  & {$72$}  & {$355$} \\
{SiO($\sigma=0.75$)} & {$4$}  & {$116.3$}  & {$12(19)$} & {$2(33)$}  & {$1(52)$}  & {$1(52)$} & {$1(3)$} & {$26.3$}  & {$85$}  & {$1,670$} \\
\end{tabular}
\label{detailed_iters_ildl}
\end{sidewaystable}

\begin{sidewaystable}[htbp]
\centering
\caption{Total number of iterations using $ILDL^{T}$ (MINRES uses the modified spd preconditioner)}
\begin{tabular}{l |c| c c c c | c | c | c }
& MINRES-CG & \multicolumn{4}{|c|}{GMRES($m$)} & FGMRES($120$) & BiCGStab & MINRES \\
Name &  & $m=20$ & $m=40$ & $m=60$  & $m=120$ & -GMRES($120$) &  &  \\ \hline
Bcsstm$10$  & $42$ & $10$ & $10$ & $10$ & $10$ & $22$  & $8.5$ & $37$ \\
Bcsstm$27$ & $8,500$ &  $\dagger$ & $\dagger$ &  $\dagger$ & $\dagger$ & $\dagger$ &  $\dagger$ & $\ddagger$ \\
Nasa$1824$ & $78$ & $31$ & $23$  & $23$ & $23$ & $74$ & $32$  & $48$\\
%iH$4$ & $125$ & $82$ & $56$ & $45$ & $45$ & $94$ & $57.0$ \\
%iNa & $105$ & $121$ & $56$ & $42$ & $42$ & $64$ & $50.5$ \\
%a$5$ & $114$ & $221$ & $61$ & $51$ & $51$ & $76$ &   $45.5$ \\
Meg$4$ & $3$ & $1$ & $1$ & $1$ & $1$ & $1$ &  $0.5$ & $2$ \\
Benzene & $23$ & $9$ & $9$ & $9$ & $9$ & $14$ &  $4.5$ & $12$ \\
Si$10$H$16$ & $288$ & $132$ & $53$  & $42$ &  $42$ & $81$ & $49.5$ & $1,192$ \\
Si$5$H$12$ & $69$ & $17$ & $17$ & $17$ & $17$ & $32$ &$11.5$ & $34$\\
SiO  & $72$ & $17$ & $17$ & $17$ & $17$ & $32$ & $11.5$ & $51$ \\
SiO ($\sigma=0.25$) & $99$ & $63$ & $31$ & $31$ & $31$ & $54$ & $40.5$ & $80$ \\
SiO ($\sigma=0.5$) & $173$ & $240$ & $42$ & $41$ & $41$  &  $85$ & $72$ & $355$ \\
SiO ($\sigma=0.75$) & $466$ & $239$ & $73$ & $52$ & $52$ & $79$ & $85$ & $1,670$
\end{tabular}
\label{total_iters_ildl}
\end{sidewaystable}

\section{Conclusions}

A two-level nested iterative scheme is proposed for solving sparse linear systems of equations where the coefficient matrix is symmetric indefinite with few negative eigenvalues. The first level is MINRES preconditioned via CG. The inner level CG is preconditioned via the original indefinite coefficient matrix. The robustness of the proposed scheme is illustrated for linear systems that arise in disk brake squeal  as well as systems that arise in a variety of test cases from the SuiteSparse Matrix Collection.

\bibliographystyle{siamplain}
\bibliography{mybib}

\begin{thebibliography}{10}

\bibitem{AxeNA14}
{\sc O.~Axelsson, M.~Neytcheva, and B.~Ahmad}, {\em A comparison of iterative
  methods to solve complex valued linear algebraic systems}, Numerical
  Algorithms, 66 (2014), pp.~811--841,
  \url{https://doi.org/10.1007/s11075-013-9764-1}.

\bibitem{BenGL05}
{\sc M.~Benzi, G.~Golub, and J.~Liesen}, {\em Numerical solution of saddle
  point problems}, Acta Numerica, 14 (2005), pp.~1--137.

\bibitem{benzi2000preconditioning}
{\sc M.~Benzi, J.~C. Haws, and M.~Tuma}, {\em Preconditioning highly indefinite
  and nonsymmetric matrices}, SIAM Journal on Scientific Computing, 22 (2000),
  pp.~1333--1353.

\bibitem{benzi1996sparse}
{\sc M.~Benzi, C.~D. Meyer, and M.~Tuma}, {\em A sparse approximate inverse
  preconditioner for the conjugate gradient method}, SIAM Journal on Scientific
  Computing, 17 (1996), pp.~1135--1149.

\bibitem{brandt1985algebraic}
{\sc A.~Brandt, S.~McCoruick, and J.~Huge}, {\em Algebraic multigrid (amg) f0r
  sparse matrix equati0ns}, Sparsity and its Applications, 257 (1985).

\bibitem{chelikowsky1994finite}
{\sc J.~R. Chelikowsky, N.~Troullier, and Y.~Saad}, {\em
  Finite-difference-pseudopotential method: Electronic structure calculations
  without a basis}, Physical review letters, 72 (1994), p.~1240.

\bibitem{choi2011minres}
{\sc S.-C.~T. Choi, C.~C. Paige, and M.~A. Saunders}, {\em Minres-qlp: A krylov
  subspace method for indefinite or singular symmetric systems}, SIAM Journal
  on Scientific Computing, 33 (2011), pp.~1810--1836.

\bibitem{DavH11}
{\sc T.~Davis and Y.~Hu}, {\em The university of florida sparse matrix
  collection}, ACM Trans. Math. Softw., 38 (2011), pp.~1:1--1:25,
  \url{https://doi.org/10.1145/2049662.2049663}.

\bibitem{duff2001algorithms}
{\sc I.~S. Duff and J.~Koster}, {\em On algorithms for permuting large entries
  to the diagonal of a sparse matrix}, SIAM Journal on Matrix Analysis and
  Applications, 22 (2001), pp.~973--996.

\bibitem{Fis11}
{\sc B.~Fischer}, {\em Polynomial Based Iteration Methods for Symmetric Linear
  Systems}, Society for Industrial and Applied Mathematics, 2011,
  \url{https://doi.org/10.1137/1.9781611971927}.

\bibitem{FoxHW48}
{\sc L.~Fox, H.~D. Huskey, and J.~H. Wilkinson}, {\em Notes on the solution of
  algebraic linear simultaneous equations}, The Quarterly Journal of Mechanics
  and Applied Mathematics, 1 (1948), pp.~149--173,
  \url{https://doi.org/10.1093/qjmam/1.1.149}.

\bibitem{freund1991polynomial}
{\sc R.~Freund}, {\em On polynomial preconditioning and asymptotic convergence
  factors for indefinite hermitian matrices}, Linear algebra and its
  applications, 154 (1991), pp.~259--288.

\bibitem{gaul2013framework}
{\sc A.~Gaul, M.~H. Gutknecht, J.~Liesen, and R.~Nabben}, {\em A framework for
  deflated and augmented {K}rylov subspace methods}, SIAM Journal on Matrix
  Analysis and Applications, 34 (2013), pp.~495--518.

\bibitem{GilMPS92}
{\sc P.~Gill, W.~Murray, D.~B. Ponceleon, and M.~A. Saunders}, {\em
  Preconditioners for indefinite systems arising in optimization}, SIAM Journal
  on Matrix Analysis and Applications, 13 (1992), pp.~292--311,
  \url{https://doi.org/10.1137/0613022}.

\bibitem{GraMQSW16}
{\sc N.~Gr{\"a}bner, V.~Mehrmann, S.~Quraishi, C.~Schr{\"o}der, and U.~{von
  W}agner}, {\em Numerical methods for parametric model reduction in the
  simulation of disk brake squeal}, ZAMM - Zeitschrift f\"ur Angewandte
  Mathematik und Mechanik, 96 (2016), pp.~1388--1405.

\bibitem{GreHL17}
{\sc C.~Greif, S.~He, and P.~Liu}, {\em Sym-ildl: Incomplete ldlt factorization
  of symmetric indefinite and skew-symmetric matrices}, ACM Trans. Math.
  Softw., 44 (2017), pp.~1:1--1:21, \url{https://doi.org/10.1145/3054948}.

\bibitem{grote1997parallel}
{\sc M.~J. Grote and T.~Huckle}, {\em Parallel preconditioning with sparse
  approximate inverses}, SIAM Journal on Scientific Computing, 18 (1997),
  pp.~838--853.

\bibitem{LehSY98}
{\sc R.~B. Lehoucq, D.~C. Sorensen, and C.~Yang}, {\em ARPACK Users Guide:
  Solution of Large-Scale Eigenvalue Problems with Implicitly Restarted Arnoldi
  Methods}, SIAM, 1998.

\bibitem{manguoglu2010weighted}
{\sc M.~Manguoglu, M.~Koyut{\"u}rk, A.~H. Sameh, and A.~Grama}, {\em Weighted
  matrix ordering and parallel banded preconditioners for iterative linear
  system solvers}, SIAM journal on Scientific Computing, 32 (2010),
  pp.~1201--1216.

\bibitem{PaiS75}
{\sc C.~C. Paige and M.~A. Saunders}, {\em Solution of sparse indefinite
  systems of linear equations}, SIAM Journal on Numerical Analysis, 12 (1975),
  pp.~617--629, \url{https://doi.org/10.1137/0712047}.

\bibitem{polizzi2009density}
{\sc E.~Polizzi}, {\em Density-matrix-based algorithm for solving eigenvalue
  problems}, Physical Review B, 79 (2009), p.~115112.

\bibitem{RozS02}
{\sc M.~Rozlozn\`ik and V.~Simoncini}, {\em Krylov subspace methods for saddle
  point problems with indefinite preconditioning}, SIAM Journal on Matrix
  Analysis and Applications, 24 (2002), pp.~368--391,
  \url{https://doi.org/10.1137/S0895479800375540}.

\bibitem{saad1993flexible}
{\sc Y.~Saad}, {\em A flexible inner-outer preconditioned gmres algorithm},
  SIAM Journal on Scientific Computing, 14 (1993), pp.~461--469.

\bibitem{Saa03}
{\sc Y.~Saad}, {\em Iterative Methods for Sparse Linear Systems}, Society for
  Industrial and Applied Mathematics, second~ed., 2003,
  \url{https://doi.org/10.1137/1.9780898718003}.

\bibitem{SaaS86}
{\sc Y.~Saad and M.~Schultz}, {\em {GMRES}: A generalized minimal residual
  algorithm for solving nonsymmetric linear systems}, SIAM Journal on
  scientific and statistical computing, 7 (1986), pp.~856--869.

\bibitem{sakurai2007cirr}
{\sc T.~Sakurai, H.~Tadano, et~al.}, {\em Cirr: a rayleigh-ritz type method
  with contour integral for generalized eigenvalue problems}, Hokkaido
  mathematical journal, 36 (2007), pp.~745--757.

\bibitem{van92}
{\sc H.~A. van~der Vorst}, {\em {Bi-CGSTAB}: A fast and smoothly converging
  variant of bi-cg for the solution of nonsymmetric linear systems}, SIAM
  Journal on Scientific and Statistical Computing, 13 (1992), pp.~631--644,
  \url{https://doi.org/10.1137/0913035}.

\bibitem{VerK13}
{\sc E.~Vecharynski and A.~Knyazev}, {\em Absolute value preconditioning for
  symmetric indefinite linear systems}, SIAM Journal on Scientific Computing,
  35 (2013), pp.~A696--A718, \url{https://doi.org/10.1137/120886686}.

\bibitem{wang2007large}
{\sc S.~Wang, E.~d. Sturler, and G.~H. Paulino}, {\em Large-scale topology
  optimization using preconditioned krylov subspace methods with recycling},
  International journal for numerical methods in engineering, 69 (2007),
  pp.~2441--2468.

\bibitem{yang2002boomeramg}
{\sc U.~M. Yang et~al.}, {\em Boomeramg: a parallel algebraic multigrid solver
  and preconditioner}, Applied Numerical Mathematics, 41 (2002), pp.~155--177.

\end{thebibliography}

\end{document}